# Dublin Descriptors

**S N Masaev[1,2,4], G A Dorrer[1,3], V V Vingert[1], E A Yakimova[1] and S V Klochkov[1]**

[1] Siberian Federal University, pr. Svobodnyj, 79, Krasnoyarsk, 660041, Russia  
[2] Control Systems LLC, 86 Pavlova Street, Krasnoyarsk, 660122, Russia  
[3] Department of System Analysis and Operations Research, Office L-409, 410, 31, Krasnoyarsky Rabochy Av., Reshetnev Siberian State University of Science and Technology, Krasnoyarsk, 660037, Russia

[4] E-mail: faberi@list.ru

**Abstract.** Dublin descriptors are under consideration. It is part of one of the global integration processes between European countries and Russia, which began in 1999. It causes a lot of controversy and approval from different sides. For the sake of clarity, an assessment is being made of the industrial application of the Dublin Descriptors. The assessment is based on the method of integral indicators. To use the method, the enterprise is formalized as a model of events at each moment in time. Each event in the enterprise is tied to the student's skill. Accordingly, students' skills are grouped by educational level. Education levels are given as Dublin descriptors. The chosen approach makes it possible to determine the correlation between levels of education and skills. It makes it possible to analyze meaningful interconnection. A universal assessment of the use of Dublin descriptors in the enterprise allows the formation of up-to-date lists of employees for higher education, professional development and training.

## 1. Introduction
Dublin Descriptors are a system of qualifications frameworks for evaluating students for bachelor's, master's, doctoral degrees. [1].

Dublin Descriptors were first proposed in 2002 as uniform requirements for higher education under the Bologna Process [1-5].

In 1999 Russia joins the Bologna Declaration and introduces three levels of higher education. It is universal descriptors, competencies and credits.

In 2003 Russia joins the development of a common framework for qualifications of the European area. It means that Russian higher education goes over three cycles. The first cycle opens the way to education in the second cycle. Accordingly, the completion of the second cycle opens access to doctoral studies. In the same year, Russia becomes a participant in the Bologna process and a member of the club. She receives the right to develop the structure of qualifications in the European higher education system.

The Bologna Process is the largest integration process on the territory of Russia and European countries. It should have end of 2010 [4]. We can only speculate what could have thwarted these plans. For example, the 2008 global financial crisis. Maybe something else.

The reasons that gave rise to the process of creating a common qualifications framework are of interest. The first reason is the different requirements of enterprises for the qualifications of graduates of educational institutions. The second reason is a huge number of educational organizations with a non-repeating system of qualifications, disciplines and credit units. The third reason is the limited mobility of university students and graduates. Lack of opportunity to read the received education for the transition to the next level of education in another country.

When will the Bologna Process be over? Will it be over? How far will European economic sanctions against Russia push him back? It is unknown. However, this process continues in Russia. Then it is relevant to consider the implementation of Dublin descriptors in Russia and evaluate their usefulness. The purpose of the study is: to assess the impact of Dublin descriptors on enterprise operations.

We going to do a tasks for that:

- Create an enterprise model;
- Link Enterprise Model to Dublin Descriptors;
- Assess the modes of operation of the enterprise (in normal mode and in the mode with Dublin descriptors);
- Analysis of results.

Enterprise model in the classical representation [6,7].

## 2. Method

Step 1. Enterprise events are set $X$. $t$ is time period is analyzed. Then a enterprise model is formed $S=\{T,X\}$. where $T=\{t:t=1,...,T_{max}\}$ - a lot of time points, $x(t)=[x^1(t),x^2(t),...,x^n(t)]^T \in X - n$ - vector of events. Events $x^i(t)$ - the value of financial expenses and income of the enterprise. The resulting system has sufficient representation $S=\{T,X\}$. You can perform standard control operations with it [8].

Step 2. Comparison of competencies of Dublin Descriptors $v(t)$ with events $x(t)$ of the enterprise $S=\{T,X\}$. We form $v_i^j$ (compliance $x_j^i$ is $v_i^j$ set as 1-yes, 0-no) from $i$ – from competence of Dublin Descriptors and $j$ event model. We get competencies of Dublin Descriptors set $V$ where $v(t) = \left[v_1^1(t), v_2^j(t), ..., v_m^n(t)\right]^T \in V$.

Payment competencies of the model is limited by resources $C$, then $C(V) \leq C$. It restriction applies to all subsystems of the researched system.

Step 3. Calculation of the integral $V$ index through the correlation matrix $R_i(x)$

$$V_i(t) = R_i(t) = \sum_{j=1}^{n} |r_{ij}(t)|, \qquad (1)$$

$$R_k(t) = \frac{1}{k-1} \overset{o}{V}_k^T(t) \overset{o}{V}_k(t) = \left\| r_{ij}(t) \right\|, \qquad (2)$$

$$r_{ij}(t) = \frac{1}{k-1} \sum_{l=1}^{k} \overset{o}{v}^i(t-l) \overset{o}{v}^j(t-l), \quad i,j=1,...,n, \qquad (3)$$

$$V_k(t) = \begin{bmatrix} v^T(t-1) \\ v^T(t-2) \\ ... \\ v^T(t-k) \end{bmatrix} = \begin{bmatrix} v^1(t-k) & v^2(t-k) & \cdots & v^n(t-k) \\ v^1(t-k) & v^2(t-k) & \cdots & v^n(t-k) \\ ... & ... & ... & ... \\ v^1(t-k) & v^2(t-k) & \cdots & v^n(t-k) \end{bmatrix}, \qquad (4)$$

where $t$ are the time instants, $r_{ij}(t)$ are the correlation coefficients of the variables $v^i(t)$ и $v^j(t)$ at the time instant $t$.

Step 4. The analysis of experimental data is performed graphically. The dynamics of the integral indicator is calculated for all periods of time

$$V = \sum_{t=1}^{T=\max} \sum_{i=1}^{n} V_i(t).$$

(5)

## 3. Characteristics of the research objects

Dublin Descriptors is being implemented to track the level of competence of enterprise personnel. The company employs 650 people. These are administrative staff and workers. The company is engaged in the processing of round timber. Has three key business processes. The first one is logging. Second. Delivery along the Yenisei River to the plant. The third is the production of round wood products. A total of 1.2 million processes have been modeled.

Table 1 shows the Dublin Descriptors implemented in the enterprise

**Table 1.** Dublin Descriptors.

| Level | Skill | Request |
|---|---|---|
| 1 Bachelor | 1.1 Knowledge and understanding | 1.1.1 Corresponds to the level of advanced textbooks, and also includes some aspects formed by knowledge of advanced positions in the field of study [1] |
| 2 Master | 2.1 Knowledge and understanding | 2.1.1 Provides a basis or opportunity for originality in the development or application of ideas, often in the context of research [1] |
| 3 PhD | 3.1 Knowledge and understanding | 3.1.1 Includes a systematic understanding of their field of study and proficiency in research methods related to this field [1] |
| 1 Bachelor | 1.2 Application of knowledge and understanding [1] | 1.2.1 By advancing and defending arguments |
| 2 Master | 2.2 Application of knowledge and understanding [1] | 2.2.1 Through the ability to solve problems, applied in a new or unfamiliar environment in a wide (or interdisciplinary) context [1] |
| 3 PhD | 3.2 Application of knowledge and understanding | 3.2.1 Demonstrates the ability to conceive, plan, implement and apply a serious research process with scientific validity. It is in the context of a contribution of original research in new areas of knowledge, carried out through the conduct of large-scale research work, the materials of which are published or referred to in national or international sources [1] |
| 1 Bachelor | 1.3. Forming judgments | 1.3.1 Includes collection and interpretation of relevant data |
| 2 Master | 2.3. Forming judgments | 2.3.1 Demonstrates ability to integrate knowledge and cope with complexities, to make judgments based on incomplete data |
| 3 PhD | 3.3. Forming judgments | 3.3.1 Requires the ability to critically analyze, evaluate and synthesize new and complex ideas [1] |
| 1 Bachelor | 1.4 Communication | 1.4.1 Communication of information, ideas, problems and solutions [1] |
| 2 Master | 2.4 Communication [1] | 2.4.1 Communicate findings and their underlying knowledge and considerations (limited coverage) to a professional and non-specialist audience (monologue) |
| 3 PhD | 3.4 Communication [1] | 3.4.1 Communication with colleagues, the wider academic community and society as a whole (dialogue) on topics related to the field of professional knowledge (wide coverage) [1] |
| 1 Bachelor | 1.5 Learning skills | 1.5.1 The skills necessary to carry out further training with a high degree of independence have been developed |
| 2 Master | 2.5 Learning skills | 2.5.1 Allow further learning to be carried out with a high degree of independence and self-regulation |
| 3 PhD | 3.5 Learning skills | 3.5.1 Will be able to contribute - in the scientific and professional context - to technological, social and cultural progress |

Simulation is performed in the software package [9].

## 4. Experiment result

Given: $n$=1.2 million values, $X$=5,641,442 thousand rubles, control is set through Dublin Descriptors ($V_{D\_descr}$). From the 7th period, one HR manager is hired in accordance with the lines of business of the enterprise to maintain the selected control syntax on it. The calculation algorithm is 447 minutes [9].

Table 2 shows the experiment result of estimating the control mode $V_i(t)$ through competence of Dublin Descriptors.

**Table 2.** Regimes: $V_{(basic\ mode)}$ and $V_{(D\_descr)}$.

| t | $V_{(basic\_mode)}$ | $V_{(D\_descr)}$ | $\Delta V$ | t | $V_{(basic\_mode)}$ | $V_{(D\_descr)}$ | $\Delta V$ |
|---|---|---|---|---|---|---|---|
| 1 | 87.34 | 110.64 | 23.30 | 30 | 96.32 | 95.32 | -1.00 |
| 2 | 70.94 | 90.26 | 19.32 | 31 | 105.10 | 104.10 | -1.00 |
| 3 | 51.43 | 69.02 | 17.59 | 32 | 98.66 | 97.66 | -1.00 |
| 4 | 56.35 | 81.32 | 24.97 | 33 | 82.19 | 81.19 | -1.00 |
| 5 | 59.26 | 89.99 | 30.73 | 34 | 76.23 | 75.22 | -1.00 |
| 6 | 73.39 | 110.23 | 36.84 | 35 | 68.52 | 68.52 | 0.00 |
| 7 | 95.25 | 123.53 | 28.28 | 36 | 60.51 | 60.51 | 0.00 |
| 8 | 92.64 | 124.07 | 31.43 | 37 | 53.13 | 53.13 | 0.00 |
| 9 | 95.53 | 125.43 | 29.90 | 38 | 61.65 | 61.65 | 0.00 |
| 10 | 70.17 | 100.32 | 30.15 | 39 | 53.51 | 53.51 | 0.00 |
| 11 | 58.42 | 74.76 | 16.35 | 40 | 51.84 | 51.84 | 0.00 |
| 12 | 56.48 | 72.97 | 16.50 | 41 | 72.03 | 72.03 | 0.00 |
| 13 | 61.88 | 79.69 | 17.81 | 42 | 93.08 | 93.08 | 0.00 |
| 14 | 71.87 | 90.79 | 18.92 | 43 | 99.23 | 99.23 | 0.00 |
| 15 | 52.45 | 67.93 | 15.48 | 44 | 115.79 | 115.78 | -0.01 |
| 16 | 53.92 | 70.08 | 16.16 | 45 | 110.12 | 110.11 | -0.01 |
| 17 | 84.06 | 101.54 | 17.48 | 46 | 103.64 | 103.64 | -0.01 |
| 18 | 114.43 | 132.62 | 18.19 | 47 | 88.22 | 87.22 | -1.00 |
| 19 | 132.20 | 132.20 | 0.00 | 48 | 69.27 | 68.27 | -1.00 |
| 20 | 153.90 | 153.92 | 0.01 | 49 | 55.14 | 54.14 | -1.00 |
| 21 | 164.54 | 164.53 | -0.01 | 50 | 63.51 | 62.51 | -1.00 |
| 22 | 150.02 | 151.00 | 0.98 | 51 | 50.02 | 49.02 | -1.00 |
| 23 | 140.74 | 144.74 | 4.00 | 52 | 61.58 | 60.58 | -1.00 |
| 24 | 115.09 | 119.09 | 3.99 | 53 | 60.33 | 59.33 | -1.00 |
| 25 | 87.02 | 91.02 | 4.00 | 54 | 147.33 | 147.33 | 0.00 |
| 26 | 100.60 | 104.60 | 4.01 | 55 | 158.41 | 158.41 | 0.00 |
| 27 | 87.59 | 91.60 | 4.00 | 56 | 156.87 | 156.87 | 0.00 |
| 28 | 76.04 | 79.04 | 3.00 | 57 | 167.90 | 167.90 | 0.00 |
| 29 | 76.26 | 76.26 | 0.00 | Total | 5 069.93 | 5 491.28 | 421.35 |

Figure 1 shows the experiment result of estimating the control mode $V_i(t)$ through Dublin Descriptors.

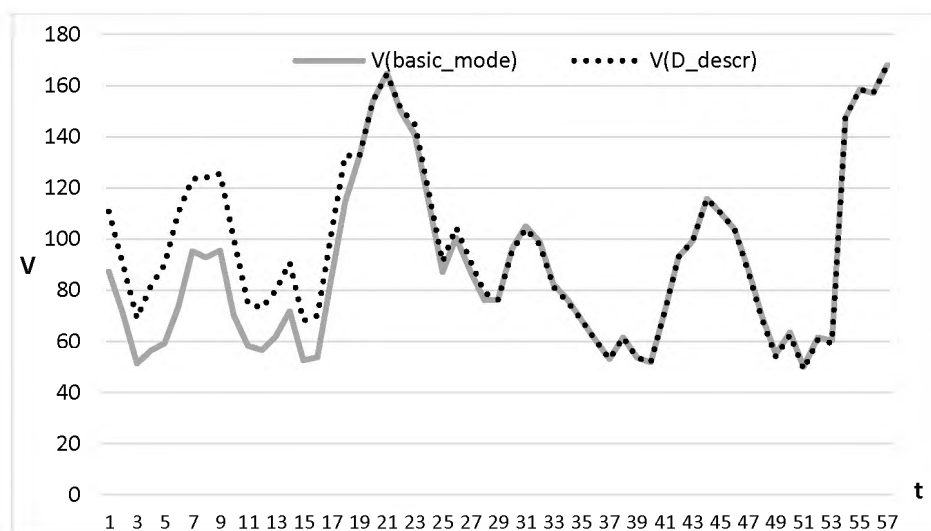

**Figure 1.** Indicator dynamics $V_i(t)$.

## 5. The discussion of the results

The installation of Dublin Descriptors for assessing the competence of personnel in the enterprise costs 28,208 thousand rubles. Then the total costs of the enterprise for five years will amount to 5,669,650 thousand rubles.

Evaluation by a universal (integral) indicator makes it possible to find competencies in demand in the industry. Determine how competencies are interconnected in practice. Describe how educational levels affect each other. It makes possible to effectively distribute the training program between universities in European countries and Russia.

## 6. Conclusion

The universal (integral) assessment indirectly reduces the time required for the implementation of the Bologna Process in Russia and Europe. It also allows you to significantly reduce the amount of processed information. Most importantly, a unified form of education will raise the level of student mobility in Russia and Europe. It is also possible to select the best educational program for training industrial workers. It makes it possible to form lists for professional development and retraining. The topical issue of the demand for Dublin descriptors in industry is considered.

The research tasks were completed:

- Created a model of the enterprise $S=\{T,X\}$;
- Associate Enterprise Model with Dublin Descriptors $v(t) = \left[v_1^1(t), v_2^j(t), ..., v_m^n(t)\right]^T \in V$;
- Assess the operating modes of the enterprise (in normal mode and in the mode with Dublin descriptors) $V_{(D\_descr)}$;
- Results analysis performed $\Delta V = V_{(D\_descr)} - V_{(basic\_mode)} = 5{,}491.28 - 5{,}069.93 = 421.35$.

The purpose of the research has been achieved.